\documentclass[10pt,twoside]{article}
\usepackage{amsmath,amsfonts,amssymb,amscd}

\font\titel=cmssbx10 scaled \magstep2
\font\naam=cmssbx10 scaled \magstep1
\font\kop=cmss10 scaled \magstep1

\newtheorem{theorem}{Theorem}
\newtheorem{corollary}{Corollary}
\newtheorem{lemma}{Lemma}

\newcommand{\CC}{{\mathbb{C}}}

\newcommand{\QQ}{{\mathbb{Q}}}

\newcommand{\HH}{{\mathbb{H}}}
\newcommand{\RR}{{\mathbb{R}}}

\newcommand{\calI}{{\cal I}}

\newcommand{\calL}{{\cal L}}

\newcommand{\calK}{{\cal K}}
\newcommand{\calM}{{\cal M}}

\begin{document}
\thispagestyle{empty}

\vspace{1cm}
\centerline{\kop DEPARTMENT OF MATHEMATICS}
\medskip
\centerline{\kop UNIVERSITY OF NIJMEGEN The Netherlands}

\vspace{3cm}

\centerline{\titel Cup products }
\medskip
\centerline{\titel and mixed Hodge structures}
\medskip
\medskip

\vspace{1.3cm}

\centerline{\naam J.H.M.~Steenbrink}

\vspace{5.3cm}
\centerline{\naam Report No. 9907 (February 1999)} 

\vspace{2.6cm}
{\kop
DEPARTMENT OF MATHEMATICS

UNIVERSITY OF NIJMEGEN

Toernooiveld

6525 ED Nijmegen

The Netherlands}

\newpage

\setcounter{page}{1}

\title{Cup products and mixed Hodge structures}
\author{J.H.M.~Steenbrink}
\date{}  
\maketitle

\begin{abstract}
Let $X$ be a complete complex algebraic variety of dimension $n$
and let $D$ be a divisor with strict normal crossings on $X$. In this paper
we show that
the cup product maps  
$$H^i(X \setminus D)\otimes H^j(X,D) \to H^{i+j}(X,D)$$ and 
$$H^i_D(X) \otimes H^j(D) \to H^{i+j}_D(X)$$ are morphisms of mixed
Hodge structures.

\end{abstract}

\maketitle
\section{Introduction}
Let $X$ and $Y$ be topological spaces, $A$ a commutative ring, $\calK^\bullet$
a complex of sheaves of $A$-modules on $X$ and $\calL^\bullet$
a complex of sheaves of $A$-modules on $Y$. Then one has the complex of
sheaves of $A$-modules $\calK^\bullet \boxtimes \calL^\bullet$ on $X \times Y$
and canonical homomorphisms
$$\HH^p(X,\calK^\bullet) \otimes_A \HH^q(Y,\calL^\bullet) \to
\HH^{p+q}(X \times Y, \calK^\bullet \boxtimes \calL^\bullet) \ .$$

In the case $X=Y$, suppose that moreover we are given a homomorphism of
complexes of sheaves of $A$-modules
$$\calK^\bullet \otimes \calL^\bullet \to \calM^\bullet$$
then, using that the restriction of $\calK^\bullet \boxtimes \calL^\bullet$ to
the diagonal of $X \times X$ is equal to $\calK^\bullet \otimes \calL^\bullet$,
one obtains homomorphisms
$$ \HH^p(X,\calK^\bullet) \otimes_A \HH^q(X,\calL^\bullet) \to
\HH^{p+q}(X,\calM^\bullet) \ .$$

As an example, we let $X$ be a smooth complete complex algebraic variety,
and consider the map $$
\mu: \Omega_X^\bullet \otimes_{\CC} \Omega_X^\bullet \to
\Omega_X^\bullet $$ given by cup product of forms. It induces the cup product
on cohomology with coefficients in $\CC$ and we can deduce Poincar\'e duality
from Serre duality in the following way. We let $n = \dim X$ and focus on
the cup product
$$
\HH^k(X,\Omega_X^\bullet) \otimes_{\CC} \HH^{2n-k}(X, \Omega_X^\bullet) \to
\HH^{2n}(X ,\Omega_X^\bullet) \simeq \CC \ .
$$
Note that the morphism $\mu$ maps $\sigma^{\geq p}\Omega_X^\bullet \otimes
\sigma^{\geq q}\Omega_X^\bullet$ to $\sigma^{\geq p+q}\Omega_X^\bullet$ so
the induced map 
$$F^p\HH^k(X,\Omega_X^\bullet) \otimes_{\CC} F^q\HH^{2n-k}(X, \Omega_X^\bullet)
\to \HH^{2n}(X, \Omega_X^\bullet)$$ is the zero map as soon as $p+q>n$. For
$p+q=n$ the pairing
$$Gr_F^p\HH^k(X,\Omega_X^\bullet) \otimes_{\CC}
Gr_F^q\HH^{2n-k}(X, \Omega_X^\bullet) \to \HH^{2n}(X ,\Omega_X^\bullet)$$
is nonsingular by Serre duality. This implies that the pairing 
$$
\HH^k(X,\Omega_X^\bullet) \otimes_{\CC} \HH^{2n-k}(X, \Omega_X^\bullet) \to
\HH^{2n}(X ,\Omega_X^\bullet) \simeq \CC \ .
$$
is also nonsingular.

Let $X$ be a complete complex algebraic variety of dimension $n$ 
and let $D$ be a divisor with strict normal crossings on $X$. In this paper
we show that
the cup product maps  
$$H^i(X \setminus D)\otimes H^j(X,D) \to H^{i+j}(X,D)$$ and 
$$H^i_D(X) \otimes H^j(D) \to H^{i+j}_D(X)$$ (the "extraordinary cup
product, cf. \cite[p. 127]{Iv} are morphisms of mixed
Hodge structures. 
A consequence of these facts is Fujiki's theorem \cite{Fu}: if $D$ is a strict 
normal crossing divisor on the smooth complete variety $X$ and $n =
\dim(X)$, then  the cup products
$$H^i(X\setminus D) \otimes H^{2n-i}(X,D) \to H^{2n}(X,D) \simeq \QQ(-n)$$
and 
$$H^i_D(X)\otimes H^{2n-i}(D) \to H^{2n}_D(X) \simeq \QQ(-n)$$
induce dualities of mixed Hodge structures, and the exact sequences of mixed
Hodge structures 
$$
\cdots \to H^i(X,D) \to H^i(X) \to H^i(D) \to H^{i+1}(X,D) \to \cdots
$$
$$ \cdots \leftarrow H^{2n-i}(X\setminus D) \leftarrow H^{2n-i}(X)
\leftarrow H^{2n-i}_D(X) \leftarrow H^{2n-i-1}(X \setminus D) \leftarrow
\cdots
$$
are dual to each other (with respect to $\QQ(-n)$). 

\section{Some cohomological mixed Hodge complexes}
A {\em cohomological mixed Hodge complex} on a variety $X$ consists of data 
\begin{itemize}
\item $(K_{\QQ}^\bullet,W)$, a complex of 
sheaves of $\QQ$-vector spaces on $X$ with
an increasing filtration $W$;
\item $(K_{\CC}^\bullet,W,F)$, a complex of sheaves of $\CC$-vector spaces on $X$
with an increasing filtration $W$ and a decreasing filtration $F$;
\item an isomorphism $\alpha: (K_{\QQ}^\bullet,W) \otimes \CC \simeq
(K_{\CC}^\bullet,W)$ in
the filtered derived category $DF(X,\CC)$ of sheaves of $\CC$-vector spaces
on $X$ with an increasing filtration
\end{itemize}
which satisfy certain axioms, see \cite[Sect.~8]{H3}. 
Such an object gives rise to mixed $\QQ$-Hodge structures on the
hypercohomology groups $\HH^k(X,K_{\CC}^\bullet)$: the $\QQ$-structure is induced by
the map defined by $\alpha$ on hypercohomology, and the weight and Hodge
filtrations are induced by $W$ and $F$ respectively. 

Let $X$ be a smooth complete variety of dimension $n$ and let $D =
\bigcup_{\alpha \in A} D_\alpha$ be a
divisor with strict normal crossings on $X$. Let $U = X \setminus D$.
Then one has a mixed Hodge structure on $H^k(U)$ defined by the "standard"
cohomological mixed Hodge complex $K^\bullet(X \log D)$ whose
$\CC$-component satisfies 
$$K^\bullet(X \log D)_{\CC} = \Omega^\bullet_X(\log D)$$
with $W$ and $F$ as in \cite[Sect.~3.1]{H2}. 

Similarly, one has a mixed Hodge structure on $H^k(D)$ defined by a standard
cohomological mixed Hodge complex $K^\bullet(D)$ which is described as
follows. 

Let $D_\bullet$ be the semisimplicial variety with $D_m$ 
the union of $(m+1)$-fold intersections of components of $D$.
We have an augmentation $\pi_\bullet:D_\bullet \to D$ which is of cohomological
descent: for any complex of sheaves of abelian groups $F^\bullet$ on $D$ one
has a quasi-isomorphism $F^\bullet \to \pi_\ast\pi^\ast F^\bullet$. (Here
$\pi_\ast\pi^\ast F^\bullet$ is the single complex associated to the double
complex $ \bigoplus_{p,q} (\pi_p)_\ast (\pi_p)^\ast F^q$). 
Then $$K^\bullet(D)_{\QQ} = \pi_\ast\pi^\ast \QQ_D = \bigoplus_p(\pi_p)_\ast
\QQ_{D_p}$$ with $$W_kK^\bullet(D)_{\QQ} = \bigoplus_{p\geq -k}(\pi_p)_\ast
\QQ_{D_p}$$ and $$K^\bullet(D)_{\CC}= \bigoplus_p(\pi_p)_\ast
\Omega^\bullet_{D_p}$$ with the filtration $W$ as in the rational case and
the filtration $F$ given by $$F^qK^\bullet(D)_{\CC} = \bigoplus_p(\pi_p)_\ast 
F^q\Omega^\bullet_{D_p}\ .$$ 

Define $\tilde{\Omega}^\bullet_D = \Omega^\bullet_D \mbox{ mod torsion }$
Then $$\tilde{\Omega}^\bullet_D \simeq
\Omega^\bullet_X/\Omega^\bullet_X(\log D)(-D)$$ with the filtration $F$
induced from $\Omega^\bullet_X$, and one has a filtered quasi-isomorphism 
$$
(\tilde{\Omega}^\bullet_D ,F) \simeq (K^\bullet(D)_{\CC},F)\ .
$$
However,  there is no weight filtration on  $\tilde{\Omega}^\bullet_D $.  

One has a natural restriction map of cohomological
mixed Hodge complexes $i^\ast: K^\bullet(X) \to K^\bullet(D)$, and 
the mixed Hodge structure on $H^\ast(X,D)$ is obtained from the
cohomological mixed Hodge complex $$K^\bullet(X,D) :=
\mbox{Cone}(i^\ast)[-1]$$ (here we take the so-called {\em mixed cone}, cf.
\cite[Sect. 3.3]{E}, and $[-1]$ denotes a shift of index in the complex). 
Observe that $$(K^\bullet(X,D)_{\CC},F) \simeq \ker((\Omega^\bullet_X,F) \to
(\tilde{\Omega}^\bullet_D ,F)) \simeq (\Omega^\bullet_X(\log D)(-D),F)$$
This observation suffices to show that the cup product mapping 
$$
H^r(X \setminus D) \otimes H^s(X,D) \to H^{r+s}(X,D)
$$
is compatible with the Hodge filtrations. Indeed, cup product induces a
morphism of complexes
$$
\Omega^\bullet_X(\log D) \otimes_\CC 
\Omega^\bullet_X(\log D)(-D) \to
\Omega^\bullet_X(\log D)(-D)
$$
which is compatible with the induced Hodge filtrations. 
Moreover observe that $\Omega^n_X(\log D)(-D) \simeq \Omega^n_X$ so
$\HH^{2n}(X,\Omega^\bullet_X(\log D)(-D)) \simeq H^n(X,\Omega_X^n) \simeq
\CC$ if $X$ is connected. Now the argument runs just like in the smooth
projective case: Serre duality gives a non-degenerate pairing between
$$Gr_F^pH^k(X \setminus D) \simeq H^{k-p}(X,\Omega_X^p(\log D))$$
and
$$Gr_F^{n-p}H^{2n-k}(X,D) \simeq H^{n-k+p}(X,\Omega_X^{n-p}(\log D)(-D)) \
.$$
However, there is no weight filtration on $\Omega^\bullet_X(\log D)(-D)$, so
we cannot us the complex $\Omega^\bullet_X(\log D)(-D)$ 
to prove compatibility of the cup product with the weight
filtration. 

\section{A weak equivalence}
The cup product map $$H^i(X \setminus D,\CC) \otimes_{\CC} H^j(X,D;\CC) \to
H^{i+j}(X,D;\CC)$$
by the results of Sect. 2 is reformulated as  a map 
$$
\HH^i(X,\Omega_X^\bullet(\log D)) \otimes_{\CC}
\HH^j(X,\Omega_{X,D}^\bullet) \to \HH^{i+j}(X,\Omega_{X,D}^\bullet) \ .
$$
However, we do not dispose of a natural cup product mapping on the level of
complexes 
$$\Omega_X^\bullet(\log D) \otimes_{\CC_X} \Omega_{X,D}^\bullet \to
\Omega_{X,D}^\bullet \ .
$$
Indeed, already $\Omega_X^\bullet(\log D) \otimes_{\CC_X} \Omega_X^\bullet$ 
does not map to $\Omega_X^\bullet$ by cup product, but to
$\Omega_X^\bullet(\log D)$. In general, for an irreducible component $C$ of
$D_m$ for some $m \geq 0$ we look for a natural target for a cup product map
on $\Omega_X^\bullet(\log D) \otimes_{\CC_X} \Omega_C^\bullet$. This is
provided by the following

\begin{lemma}
Let $C$ be an irreducible component of $D_m$ for some $m$. Then $C$ is
a smooth subvariety of $X$. 
Let $\calI_C \subset O_X$ denote its ideal sheaf.
Then $\calI_C\Omega_X^\bullet(\log D)$ is a subcomplex of 
$\Omega_X^\bullet(\log D)$. 
\end{lemma}
{\em Proof.} 
Let $P \in C$. Choose local holomorphic coordinates $(z_1,\ldots,z_n)$ on
$X$ centered at $P$ such that $\calI_{C,P} = (z_1,\ldots,z_k)O_{X,P}$ and
$\calI_{D,P} = (z_1\cdots z_l)O_{X,P}$ for some $k \leq l \leq n$. For
$\omega \in \calI_C\Omega_X^p(\log D)_P$ write $\omega = \sum_{i=1}^k
z_i\omega_i$ with $\omega_i \in \Omega_X^p(\log D)_P$ for $i=1,\ldots,k$.
Then $$d\omega = \sum_{i=1}^k z_i(\frac{dz_i}{z_i} \wedge \omega_i +
d\omega_i) \in \calI_C\Omega_X^{p+1}(\log D)_P \ .$$

\vspace{4mm}
\noindent
We denote the quotient complex $\Omega_X^\bullet(\log
D)/\calI_C\Omega_X^\bullet(\log D)$ by $\Omega_X^\bullet(\log D) \otimes 
O_C$. We equip it with the filtrations $W$ and $F$ as a quotient of
$\Omega_X^\bullet(\log D)$. 

\begin{theorem}
The complex $\Omega_X^\bullet(\log D)\otimes O_C$ is quasi-isomorphic to
$i_C^\ast\RR j_\ast\CC_{X \setminus D}$ where $i_C:C \to X$ and $j:X
\setminus D \to X$ are the inclusion maps.
\end{theorem}
{\em Proof.} 
We have an isomorphism $(\Omega_X^\bullet(\log D),W) \simeq (\RR
j_\ast\CC_{X \setminus D}, \tau_\leq)$ in $DF(X,\CC)$, which by restriction
to $C$ gives an isomorphism 
$$(i_C^\ast\Omega_X^\bullet(\log D),W) \simeq
(i_C^\ast\RR j_\ast\CC_{X \setminus D}, \tau_\leq)$$ in $DF(C,\CC)$.  
It remains to be proven that the quotient map 
$$(i_C^\ast\Omega_X^\bullet(\log D),W) \to (\Omega_X^\bullet(\log
D)\otimes O_C,W)$$ is a filtered quasi-isomorphism. 
To deal with this problem, note that for all $k \geq 0$ 
one has the Poincar\'e residue map 
$$R_k:Gr^W_k \Omega_X^\bullet(\log D) \to
(\pi_{k-1})_\ast\Omega^\bullet_{D_{k-1}}[-k]$$ (where $D_{-1}:= X$) 
which is an isomorphism of complexes. It has components $$R_I: Gr^W_k
\Omega_X^\bullet(\log D) \to \Omega^\bullet_{D_I}[-k]$$ where $I$ is a
subset of $A$ of cardinality $k$ and $D_I := \bigcap_{a\in I} D_a$. It
follows that 
$$i_C^\ast Gr^W_k \Omega_X^\bullet(\log D) \simeq \bigoplus_{\sharp I =
k} i_C^\ast \CC_{D_I}[-k] \simeq \bigoplus_{\sharp I =
k} \CC_{D_I \cap C}[-k] \ .$$ 
Claim: the image of $\calI_C\Omega^p_X(\log D) \cap
W_k\Omega^p_X(\log D)$ under the map $R_I$ coincides with
$\calI\Omega^{p-k}_{D_I} + d\calI_C\wedge\Omega^{p-k-1}_{D_I}$. 

Assuming the claim, we find that $R_k$ induces an isomorphism 
$$Gr^W_k\Omega^\bullet_X(\log D) \otimes O_C \simeq \bigoplus_{\sharp I
= k} \Omega_{D_I \cap C}^\bullet[-k] \ .$$
Let us prove the claim. The map $R_I$ presupposes an ordering of the set $A$
of irreducible components of $D$. Write $I = \{i_1,\ldots,i_k\}$ with $i_1 <
\ldots < i_k$ and choose local coordinates $(z_1,\ldots,z_n)$ on $X$
centered at $P \in C$ such that $D_{i_r}$ is defined near $P$ by $z_r=0$ for
$r=1,\ldots,k$ and $\calI_{C,P}$ is generated by $z_j$ for $j \in J$. Put
$J_1 = J \cap \{1,\ldots,k\}$ and $J_2 = J \setminus J_1$. Also suppose that
$D$ is defined near $P$ by $z_1\cdots z_l= 0$. Then $l \geq k$ and $J
\subset \{1,\ldots,l\}$. 

For $j \in J_2$ choose $\eta_i \in \Omega^{p-k-1}_{D_I,P}$ and $\zeta_i \in
\Omega^{p-k}_{D_I,P}$ with lifts $ \tilde{\eta}_i$ and $\tilde{\zeta}_i$ in
$\Omega^{p-k-1}_{X,P}$ and $\Omega^{p-k}_{X,P}$ respectively. Let 
$$\omega =
\sum_{j \in J_2} \frac{dz_1}{z_1} \wedge \cdots \wedge \frac{dz_k}{z_k}
\wedge (dz_j \wedge \tilde{\eta}_i + z_j\tilde{\zeta}_j) \ .$$
Then $\omega \in \calI_C\Omega^p_X(\log D)_P \cap
W_k\Omega^p_X(\log D)_P$ and $$R_I(\omega) = \sum_{j \in J_2} (dz_j\wedge
\eta_j + z_j\zeta_j) \ .$$ Also remark that $R_{I^\prime}(\omega)= 0$ if $I
\neq I^\prime \subset A$ with $\sharp I^\prime = k$. 
Hence we have the inclusion 
$$\bigoplus_I(\calI\Omega^{p-k}_{D_I,P} +
d\calI_C\wedge\Omega^{p-k-1}_{D_I,P}) \subset 
R_k(\calI_C\Omega^p_X(\log D)_P \cap W_k\Omega^p_X(\log D)_P) \ .$$

To prove the reverse inclusion, we let $\xi_i = \frac{dz_i}{z_i}$ if $1 \leq
i \leq l$ and $\xi_i = dz_i$ if $i >l$. Also, for $B = \{b_1,\ldots,b_r\}
\subset \{1,\ldots,n\}$ with $b_1 < \cdots < b_r$ we put $\xi_B = \xi_{b_1}
\wedge \cdots \wedge \xi_{b_r}$. With this notation, $\Omega^p_X(\log D)_P$
is the free $O_{X,P}$-module with basis the $\xi_B$ with $\sharp B=p$. We
have 
$$\calI_C\Omega^p_X(\log D)_P = \bigoplus_{\sharp B = p} \calI_{C,P}\xi_B$$
and 
$$W_k\Omega^p_X(\log D)_P = \bigoplus_{\sharp B = p} W_k\Omega^p_X(\log D)_P
\cap O_{X,P}\xi_B = \bigoplus_{\sharp B = p} J(B,k)\xi_B$$
where $J(B,k)$ is an ideal of $ O_{X,P}$ generated by squarefree monomials. 

For any $B$ and any squarefree monomial $z_E \in J(B,k)$ with $R_I(z_E\xi_B)
\neq 0$ one has $\{1,\ldots,k\} \subset B$ and $B \cap \{k+1,\ldots,l\}
\subset E$. If moreover $z_E \in \calI_{C,P}$ then $J_2 \cap E \neq
\emptyset$. Choose $j \in J_2 \cap E$. If $j \in B$ then 
$$z_E\xi_B = \pm z_E \frac{dz_j}{z_j}\wedge \xi_{B \setminus \{j\}} = \pm
dz_j \wedge z_{E \setminus \{j\}} \xi_{B \setminus \{j\}} \in  
d\calI_C\wedge\Omega^{p-1}_X(\log D)_P$$ so $R_I(z_E\xi_B ) \in
d\calI_C\wedge\Omega^{p-k-1}_{D_I}$. On the other hand, if $j \not\in B$ then 
$$R_I(z_E\xi_B) = z_jR_I(z_{E \setminus \{j\}} \xi_B) \in
\calI_C\Omega^{p-k}_{D_I} \ .$$

\begin{corollary}
$\HH^k(C,\Omega_X^\bullet(\log D)\otimes O_C) \simeq H^k(U_C \setminus
D,\CC) $ where $U_C$ is a tubular neighborhood of $C$ inside $X$.
\end{corollary}

\begin{corollary}
One has a cohomological mixed Hodge complex $K^\bullet(C \log D)$ on $C$
with $(K^\bullet(C \log D)_{\QQ},W) = (i_C^\ast\RR j_\ast\QQ_{X \setminus
D}, \tau_{\leq})$ and $(K^\bullet(C \log D)_{\CC},W,F) =
(\Omega_X^\bullet(\log D)\otimes O_C,W,F)$. This defines a mixed Hodge
structure on $H^k(U_C \setminus D,\CC) $. Moreover, $W_0K^\bullet(C \log D)
\simeq K^\bullet(C)$ so $$W_kH^k(U_C \setminus D) = \mbox{Image of } [H^k(C)
\simeq H^k(U_C) \to H^k(U_C \setminus D)]\ .$$
\end{corollary}

The data of all $K^\bullet(D_I \log D)$ for $I \subset A$ give rise to a
cohomological mixed Hodge complex on the semi-simplicial variety
$D_\bullet$. We define $$K^\bullet(D \log D) = (\pi_\bullet)_\ast
K^\bullet(D_\bullet \log D) \ .$$ This is a cohomological mixed Hodge complex on
$D$ such that $K^\bullet(D \log D)_{\QQ} \simeq i_D^\ast Rj_\ast\QQ_{X
\setminus D}$. It gives a mixed Hodge structure on $H^k(U_D \setminus
D)$ where $U_D$ is a tubular neighborhood of $D$. The spectral sequence 
$$E_1^{pq} = \HH^q(D_p,K^\bullet(D_p \log D)) \Rightarrow
\HH^{p+q}(D,K^\bullet(D \log D))$$
can be considered as the Mayer-Vietoris spectral sequence corresponding to a
covering of $U_D \setminus D$ by deleted neighborhoods $U_{D_i} \setminus
D$. Observe that we dispose of a natural morphism of cohomological mixed
Hodge complexes $K^\bullet(X \log D) \to K^\bullet (D \log D)$ which on
cohomology induces the restriction mapping $H^k(X \setminus D) \to H^k(U_D
\setminus D)$. 
We now define 
\begin{equation}
\tilde{K}^\bullet(X,D) = \mbox{cone}(K^\bullet(X \log D) \to K^\bullet(D
\log D))[-1] 
\end{equation}

Note that the inclusions $K^\bullet(X) \to K^\bullet(X \log D)$ and
$K^\bullet(D) \to K^\bullet(D \log D)$ induce a morphism of cohomological
mixed Hodge complexes 
$$\beta: K^\bullet(X,D) \to \tilde{K}^\bullet(X,D) \ .$$

\begin{lemma}
The map induced by $\beta$ on cohomology is a quasi-isomorphism.
\end{lemma}
{\em Proof.} By excision, the map $H^k(X,D) \simeq H^k(X,U_D) 
\to H^k(X \setminus D,U_D \setminus D)$ is an isomorphism for all $k$. 

\noindent {\bf Remark} \hspace{4mm} As this lemma is true also locally on
$X$, we may even conclude that $\beta$ is a quasi-isomorphism. 

\begin{corollary} The cohomological mixed Hodge complexes $K^\bullet(X,D)$ and 
$\tilde{K}^\bullet(X,D)$ determine the same mixed Hodge structure on
$H^k(X,D)$. 
\end{corollary}
Indeed, $\beta$ induces a morphism of mixed Hodge structures which is an
isomorphism of vector spaces, hence an isomorphism of mixed Hodge
structures. 

\vspace{4mm}
\noindent Now we proceed to the definition of the cup product on the level
of complexes. Write $\tilde{\Omega}^\bullet_{X,D} =
\tilde{K}^\bullet(X,D)_{\CC}$. 
For each component $C$ of $D_\bullet$ we have a natural cup
product 
$$\mu_C: \Omega_X^\bullet(\log D) \otimes_{\CC} \Omega^\bullet_C \to
\Omega_X^\bullet(log D) \otimes O_C \ .$$
These glue to give a cup product  
$$\mu: \Omega_X^\bullet(\log D) \otimes_{\CC} \Omega^\bullet_{X,D} \to
\tilde{\Omega}^\bullet_{X,D} $$
which is compatible with the filtrations $W$ and $F$. We conclude

\begin{theorem}
The cup product maps 
$$H^i(X \setminus D)\otimes H^j(X,D) \to H^{i+j}(X,D)$$ 
are morphisms of mixed Hodge structures.
\end{theorem}
{\bf Remark} \hspace{4mm} If $Y$ is a complete complex algebraic variety
with a closed subvariety $Z$ such that $Y \setminus Z$ is smooth, then there
exists a proper modification $f:X \to Y$ such that $X$ is smooth, $f$ maps
$X \setminus f^{-1}(Z)$ isomorphically to $Y \setminus Z$ and $D :=
f^{-1}(Z)$ is a divisor with strict normal crossings on $X$. Then one has
isomorphisms of mixed Hodge structures $H^i(Y \setminus Z) \to H^i(X
\setminus D)$ and $H^j(Y,Z) \to H^j(X,D)$ so that case is reduced to the
strict normal crossing case. 

\vspace{4mm}
\noindent {\bf Remark} \hspace{4mm} The restrictions of $K^\bullet(X
\log D)$, $K^\bullet(X,D)$ and $\tilde{K}^\bullet(X,D)$ are all equal to
$K^\bullet(X \setminus D)$. Now consider the following situation: $Y$ is a
complete complex algebraic variety with closed subvarieties $Z$ and $W$ such
that $Y \setminus (Z \cup W)$ is smooth and $Z \cap W = \emptyset$. 
Then there is a cup product 
$$H^i(Y \setminus Z,W) \otimes H^j(Y \setminus W,Z) \to H^{i+j}(Y,Z \cup W)$$
which is a morphism of mixed Hodge structures and induces a perfect duality
if $i+j=2\dim(Y)$. The proof uses a proper modification $f:X \to Y$ such
that $X$ is smooth, $f$ maps
$X \setminus f^{-1}(Z \cup W)$ isomorphically to $Y \setminus (Z \cup W)$
and $D = f^{-1}(Z)$ and $E = f^{-1}(W)$ are divisors with normal crossings
on $X$. By a glueing process one obtains cohomological mixed Hodge complexes
$K^\bullet(X \log D,E)$ etc. and a cup product map 
$$K^\bullet(X \log D,E)_\CC \otimes K^\bullet(X \log E,D)_\CC \to
\tilde{K}^\bullet(X,D \cup E)_\CC$$
which is compatible with $W$ and $F$. 

This answers a question raised to me by V.~Srinivas. 
\section{The extraordinary cup product}
Again, let $X$ be a complete smooth complex algebraic variety and let $D$ be a
divisor with strict normal crossings on $X$. The local cohomology groups
$H^k_D(X) = H^k(X,X \setminus D)$ get a mixed Hodge structure using the
cohomological mixed Hodge complex
\begin{equation}
K^\bullet_D(X) = \mbox{cone}[K^\bullet(X) \stackrel{u}{\to} K^\bullet(X \log D)]
\end{equation}
but by excision we may as well take 
\begin{equation}
\tilde{K}^\bullet_D(X) = \mbox{cone}[K^\bullet(D) \stackrel{v}{\to} K^\bullet(D \log D)]
\end{equation}
Observe that the morphisms $u_\CC$ and $v_\CC$ are injective, even after
taking $Gr_FGr^W$, so we have bifiltered quasi-isomorphisms 
$$(K^\bullet_D(X)_\CC,W,F) \to (\mbox{coker}(u_\CC),W,F)[-1]$$ and 
$$\tilde{K}^\bullet_D(X)_\CC,W,F) \to (\mbox{coker}(v_\CC),W,F)[-1] \ .$$
Moreover, the natural cupproduct 
$$K^\bullet(X\log D)_\CC \otimes_\CC K^\bullet(D)_\CC \to K^\bullet(D \log
D)_\CC$$ 
maps $K^\bullet(X)_\CC \otimes_\CC K^\bullet(D)_\CC$ to $K^\bullet(D)_\CC$,
so induces a cup product map 
\begin{equation}
\mbox{coker}(u_\CC) \otimes_\CC K^\bullet(D)_\CC \to \mbox{coker}(v_\CC)
\end{equation}
which is compatible with the filtrations $W$ and $F$. Hence we conclude 

\begin{theorem} The extraordinary cup product map
$$H^i_D(X) \otimes H^j(D) \to H^{i+j}_D(X)$$ 
is a morphism of mixed Hodge structures.
\end{theorem}

\end{document}